\documentclass[a4paper,11pt]{article}
\pdfoutput=1

\usepackage[a4paper,tmargin=3truecm,bmargin=3truecm,rmargin=2.5truecm,
lmargin=2.5truecm,twoside,verbose=true]{geometry}
\usepackage{cancel,graphicx}

\usepackage{amsmath,amssymb}
\usepackage[amsmath, hyperref, thmmarks]{ntheorem}
\usepackage[all]{xypic}
\usepackage[pdftex]{hyperref}
\usepackage[english]{babel}

\numberwithin{equation}{section}

\allowdisplaybreaks[1]

\def\question#1{}

\newcommand\caB{{\mathcal B}}
\newcommand\caC{{\mathcal C}}

\newcommand\caL{{\mathcal L}}

\newcommand\gone{{ \mathchoice {1\mskip-4mu\mathrm{l} } {1\mskip-4mu\mathrm{l} }{1\mskip-4.5mu\mathrm{l} } {1\mskip-5mu\mathrm{l}} }}
\newcommand\gR{{\mathbb R}}

\newcommand\gC{{\mathbb C}}

\newcommand\gB{{\mathbb B}}
\newcommand\gN{{\mathbb N}}
\newcommand\gZ{{\mathbb Z}}

\newcommand\algA{{\mathbf A}}

\newcommand\ehH{\mathcal H}

\newcommand\eps{{\varepsilon}}

\newcommand\fois{\mathord{\cdot}}

\newcommand\dd{{\text{\textup{d}}}}

\newcommand\norm{\mathord{\parallel}}
\newcommand\defin{\bf}

\theoremsymbol{}
\theorembodyfont{\slshape}
\theoremheaderfont{\normalfont\bfseries}
\theoremseparator{}
\newtheorem{Theorem}{Theorem}[section]
\newtheorem{theorem}[Theorem]{Theorem}

\newtheorem{proposition}[Theorem]{Proposition}

\theorembodyfont{\upshape}
\theoremsymbol{\ensuremath{\blacklozenge}}

\newtheorem{example}[Theorem]{Example}

\newtheorem{definition}[Theorem]{Definition}

\theoremstyle{nonumberplain}
\theoremheaderfont{\scshape}
\theorembodyfont{\normalfont}
\theoremsymbol{\ensuremath{\blacksquare}}

\qedsymbol{\ensuremath{_\blacksquare}}
\theoremclass{LaTeX}

\title{Noncommutative Supergeometry and Quantum Supergroups\footnote{Work
supported by the Belgian Interuniversity Attraction Pole (IAP) within the framework ``Dynamics, Geometry and Statistical Physics'' (DYGEST), and by the Max Planck Institut f\"ur Mathematick Bonn.}}
\date{}
\author{Axel de Goursac}

\begin{document}

\maketitle
\vspace*{-1cm}
\begin{center}
\textit{Charg\'e de Recherches au FRS-FNRS,\\ IRMP, Universit\'e Catholique de Louvain,\\ Chemin du Cyclotron 2, 1348 Louvain-la-Neuve, Belgium\\
e-mail: \texttt{axelmg@melix.net}}\\
\end{center}%

\vskip 2cm

\begin{abstract}
This is a review of concepts of noncommutative supergeometry - namely Hilbert superspace, C*-superalgebra, quantum supergroup - and corresponding results. In particular, we present applications of noncommutative supergeometry in harmonic analysis of Lie supergroups, non-formal deformation quantization of supermanifolds, quantum field theory on noncommutative spaces; and we give explicit examples as deformation of flat superspaces, noncommutative supertori, solvable topological quantum supergroups.
\end{abstract}

\vskip 1cm
%
%
%

\section*{Introduction}

Differential supergeometry has its origin \cite{Berezin:1976,Kostant:1977} in the second half of the 20th century. Its objects are supermanifolds involving, besides the usual commuting coordinates, also anticommuting coordinates. The corresponding algebras of functions are then $\gZ_2$-graded commutative. This graded point of view is very useful in various domains of mathematics such as differential geometry, spinor geometry, Clifford algebras, and also in physics for fermionic or supersymmetric theories, or for the BRST formalism in gauge theories.

Besides, noncommutative geometry is a growing field of mathematics \cite{Connes:1994} whose essential principle lies in the duality between spaces and commutative algebras, so that the properties of spaces can be characterized algebraically. Then, a noncommutative algebra can be seen as corresponding to some ``noncommutative space''. This very rich way of thinking allows generalizing classical notions and theorems of usual geometry, and it is even possible to prove some new results for differential geometry in this more general noncommutative framework. Moreover, noncommutative geometry is the appropriate framework for quantization and also a challenging tool to build new models in quantum physics. It has in particular deep links with quantum gravity, quantum Hall effect and string theories \cite{Seiberg:1999vs,Douglas:2001ba}.

Putting together these two ingredients -- supergeometry and noncommutativity -- gives rise to noncommutative supergeometry (NCSG). Noncommutative supergeometry deals with (noncommutative) $\gZ_2$-graded algebras that correspond to ``noncommutative superspaces'' in the point of view of noncommutative geometry described above. It should be the right framework for fermionic or supersymmetric quantum field theories on noncommutative spacetime, or BRST formalism for noncommutative spaces... Early results \cite{deGoursac:2010zb,deGoursac:2008bd} concerning an application of NCSG to quantum field theory on the Moyal space also motivated the author to develop this program of NCSG at a mathematical level as well as to look at physical applications.

Of course, there is already an extensive literature on graded algebras (see e.g. \cite{Nastasescu:1982}), useful for NCSG at the algebraic level. Some work in the direction of noncommutative $Q$-manifolds was also performed in \cite{Schwarz:2000,Schwarz:2003}. However, at the analytic level, NCSG was really initiated in \cite{Bieliavsky:2010su} with the definition of its basic objects: Hilbert superspaces and C*-superalgebras that endow noncommutative superspaces with topological properties. In this paper, we review these analytic notions, as well as the corresponding notion of symmetry -- quantum supergroup -- defined and studied in the Fr\'echet setting \cite{deGoursac:2014kv}. We stress the relevance of this definition of Hilbert superspace, rather than the standard one, by presenting an extended Stone--von Neumann theorem valid in any signature. At the operator algebra level, a universal deformation formula that extends Rieffel's deformation quantization to the graded setting is valid in the category of C*-superalgebras, underlying its relevance. We also look at an application of NCSG to quantum field theory.

\section{Classical supergeometry}
\label{sec-concrete}

In this section, we recall basics of supergeometry (see \cite{Berezin:1976,Kostant:1977,Leites:1980,Manin:1997,Deligne:1999su,Carmeli:2011} for the sheaf-theoretic approach and \cite{DeWitt:1984,Rogers:2007,Tuynman:2005}) for the concrete approach).

A {\defin supermanifold} $M$ of dimension $m|n$ is the data of a smooth manifold of dimension $m$ called its body and denoted $\gB M$ together with a sheaf on $\gB M$ of $\gZ_2$-graded commutative algebras; and this locally ringed space is locally isomorphic to $\caC^\infty_{\gR^m}\otimes\bigwedge\gR^n$. We denote by $\caC^\infty(M)$ the algebra of global sections of the sheaf corresponding to $\gB M$. It is a $\gZ_2$-graded commutative Fr\'echet algebra. The simplest example is the locally ringed space $\caC^\infty_{\gR^m}\otimes\bigwedge\gR^n$ itself, called {\defin superspace} of dimension $m|n$ and denoted by $\gR^{m|n}$, for which we have
\begin{equation}
\caC^\infty(\gR^{m|n})\simeq \caC^\infty(\gR^m)\otimes\bigwedge \gR^n. \label{eq-decomp}
\end{equation}

To obtain explicit expressions, let us fix the canonical basis $(\xi^i)_{1\leq i\leq n}$ of $\gR^n$. A smooth superfunction $f$ on the superspace $\gR^{m|n}$, i.e. a section $f\in\caC^\infty(\gR^{m|n})$, can be decomposed as
\begin{equation*}
\sum_{I}f_I(x)\xi^I,
\end{equation*}
for any $x\in\gR^m=\gB \gR^{m|n}$, and where the sum is taken on all ordered subsets $I$ of $\{1,\dots,n\}$ with $\xi^I:=\prod_{i\in I}\xi^i =\xi^{i_1} \xi^{i_2} \cdots \xi^{i_k}$ if $I=\{i_1, \dots, i_k\}$ (ordered), and $\xi^\emptyset = 1$ by convention; the coefficients $f_I$ being smooth functions in $\caC^\infty(\gR^m)$. 

Owing to the equivalence of this definition with the one of the concrete approach, we will note in the following that $f(x,\xi)=\sum_{I}f_I(x)\xi^I$ for any $(x,\xi)\in\gR^{m|n}$, even if a rigorous definition of this notation involves functors of points \cite{functors}. For any two (ordered) subsets $I = \{i_1, \dots, i_l\}$ and $J = \{j_1, \dots, j_\ell\}$ of $\{1,\dots,n\}$ we define $\eps(I,J)$ to be zero if $I$ and $J$ overlap; if $I\cap J=\emptyset$, we set $\eps(I,J)$ to the parity of the list $(i_1, \dots, i_k, j_1, \dots, j_\ell)$, defined as $-1$ raised to the number of transpositions needed to put it in increasing order. It satisfies
\begin{equation}
\eps(I,J)=(-1)^{|I||J|}\eps(J,I),\qquad \eps(I,J\cup K)=\eps(I,J)\eps(I,K)\text{ if }J\cap K=\emptyset.\label{eq-eps}
\end{equation}
As a consequence, we can perform explicit products of superfunctions $\xi^I\fois\xi^J=\eps(I,J)\xi^{I\cup J}$ for $I\cap J=\emptyset$. We recall the Lebesgue--Berezin {\defin integration} for superfunctions:
\begin{equation}
\int_{\gR^{m|n}}\dd x\dd\xi\, f(x,\xi):=\int_{\gR^m}\dd x\,f_{\{1,\dots,n\}}(x).\label{eq-berezin}
\end{equation}

\section{Philosophy of NCSG}

The duality between spaces and commutative algebras is transverse to numerous domains of Mathematics. At the topological level, the Gelfand--Na\"imark theorem states that there is a contravariant correspondence between the category of locally compact Hausdorff spaces and the one of commutative C*-algebras: $X\mapsto \caC(X)=\{X\to \gC\text{ continuous}\}$ defining an equivalence of categories. Then, the category of noncommutative topological spaces can be defined as the dual of the category of C*-algebras \cite{Connes:1994}.

Omitting topological conditions and in view of the $\gZ_2$-graded commutative structure of the functions $\caC^\infty(M)$ on a supermanifold $M$, we can define {\defin noncommutative superspaces} as associated to $\gZ_2$-graded (noncommutative) algebras. So the difference with noncommutative geometry is simply a $\gZ_2$-grading, but we will see that this grading induces important features. This algebraic definition of noncommutative superspaces is very simple and gives rise to a well-known object (see for example \cite{Nastasescu:1982} on graded algebras). Note that geometric objects were considered and constructed for noncommutative superspaces \cite{deGoursac:2008bd} such as the noncommutative analogues of vector fields (graded derivations), de Rham forms (bigraded noncommutative differential calculus based on graded derivations), connections and covariant derivatives (noncommutative $\eps$-connections), curvatures, gauge transformations... These tools were actually defined for graded algebras with a more general grading and a bicharacter (like color algebras), opening the subject of {\defin noncommutative graded} (or {\defin color}) {\defin geometry} \cite{deGoursac:2008bd}.

\medskip

However, topological structures like C*-algebras play a crucial role in noncommutative geometry. For NCSG, we introduced in \cite{Bieliavsky:2010su} two topological notions, that of Hilbert superspaces and of C*-superalgebra that we recall here.

The standard notion of Hilbert superspace \cite{Deligne:1999su,Carmeli:2006}, also used in \cite{Alldridge:2013} consists of the data of a $\gZ_2$-graded Hilbert space $\ehH=\ehH_0\oplus\ehH_1$ with hermitian positive definite scalar product $(-,-)$ and together with an inner product defined by
\begin{equation}
\langle \varphi,\psi\rangle:= i^{|\varphi||\psi|}(\varphi,\psi).\label{eq-innerprod2}
\end{equation}
This inner product is of degree 0 and superhermitian, i.e.~$\forall x,y\in\ehH$, $\overline{\langle x,y\rangle}=(-1)^{|x||y|}\langle y,x\rangle$. The notion of Hilbert superspace presented below is a generalization of the one associated to \eqref{eq-innerprod2}. Indeed, the inner product can be of degree 1 below and the link between $\langle-,-\rangle$ is by a more general operator $J$, not only $J(\psi)=i^{|\psi|}\psi$.
\begin{definition}[\cite{Bieliavsky:2010su}]
\label{def-hs}
A {\defin Hilbert superspace} of parity $n\in\gZ_2$ is a complex $\gZ_2$-graded vector space $\ehH=\ehH_+\oplus\ehH_-$ endowed with a superhermitian (sesquilinear) homogeneous inner product $\langle-,-\rangle$ of degree $n$, such that there exists a {\defin fundamental symmetry} $J$ of degree $n$ i.e.~an endomorphism $J$ of $\ehH$ satisfying $\forall x,y\in\ehH$,
\begin{itemize}
\item $J^2(x)=(-1)^{(n+1)|x|}x$, and $\langle J(x),J(y)\rangle=\langle x,y\rangle$,
\item $(x,y)_J:=\langle x,J(y)\rangle$ defines a hermitian positive definite scalar product on $\ehH$ for which $\ehH$ is complete.
\end{itemize}
\end{definition}

The Hilbert topology of a Hilbert superspace $\ehH$ is actually independent of the choice of the fundamental symmetry $J$ \cite{deGoursac:2014svn}, so the space of continuous endomorphisms $\caB(\ehH)$ is canonically defined. In parity 1, a Hilbert superspace $(\ehH,\langle-,-\rangle)$ is a Krein space, while a Hilbert superspace of parity 0 corresponds to the orthogonal sum of the two Krein spaces $(\ehH_+,\langle-,-\rangle)$ and $(\ehH_-,i\langle-,-\rangle)$. This category of Hilbert superspaces is stable for the direct sum and the tensor product. The next Example produces Hilbert superspaces in the sense of Definition 
\ref{def-hs} but not in the sense of \eqref{eq-innerprod2} for $n>0$.
\begin{example}[\cite{Bieliavsky:2010su,deGoursac:2014svn}]
\label{ex-l2}
In view of the decomposition \eqref{eq-decomp}, we set $L^2(\gR^{m|n})=L^2(\gR^m)\otimes \bigwedge \gR^n$. It is a Hilbert superspace of parity $n$ (mod 2) with the usual $\gZ_2$-grading and the superhermitian scalar product
\begin{equation}
\langle f_1,f_2\rangle=\int_{\gR^{m|n}}\dd x\dd\xi\,\ \overline{f_1(x,\xi)}f_2(x,\xi)
\end{equation}
by using the Berezin integration \eqref{eq-berezin}. A fundamental symmetry is given by the analog of the Hodge operation $J(f)(x,\xi)=\sum_I \eps(I,\complement I) f_I(x)\xi^{\complement I}$, where $\eps$ was defined in \eqref{eq-eps} and $\complement I$ denotes the complement of the subset $I$ in $\{1,\dots,n\}$.
\end{example}

For any bounded operator $T\in\caB(\ehH)$, there exists a {\defin superadjoint} $T^\dag\in\caB(\ehH)$ defined by
\begin{equation}
\forall x,y\in\ehH\quad:\quad \langle T^\dag(x),y\rangle=(-1)^{|T||x|}\langle x,T(y)\rangle.
\end{equation}
The next Definition corresponds to the standard notion of real form in \cite{Deligne:1999su}.
\begin{definition}[\cite{Bieliavsky:2010su}]
A {\defin superinvolution} on a complex $\gZ_2$-graded algebra $\algA$ is an antilinear map $\dag:\algA\to\algA$ homogeneous of degree 0 such that
\begin{equation*}
\forall a,b\in\algA\quad:\quad (a^\dag)^\dag=a,\qquad (ab)^\dag=(-1)^{|a||b|}b^\dag a^\dag.
\end{equation*}
For example, $\caB(\ehH)$ is a $\gZ_2$-graded algebra with a superinvolution given by the superadjoint.

A {\defin C*-superalgebra} is a complex $\gZ_2$-graded algebra with superinvolution, faithfully represented on a Hilbert superspace $\ehH$ by $\rho:\algA\to\caB(\ehH)$ (with $\rho$ of degree 0 and compatible with superinvolutions), and such that $\algA$ is closed for the operator norm topology of $\caB(\ehH)$ via $\rho$.
\end{definition}
This concrete definition allows directly to define the tensor product of two C*-superalgebras. Spectral properties of C*-superalgebras are currently under study. Morphisms of C*-superalgebras are homogeneous homomorphisms of degree 0 compatible with superinvolutions. Owing to the beginning of this section, we can define the category of {\defin noncommutative topological superspaces} as the dual category of C*-superalgebras.
\begin{example}[\cite{Bieliavsky:2010su}]
The space $L^\infty(\gR^{m|n})=L^\infty(\gR^m)\otimes\bigwedge\gR^n$ is a supercommutative C*-superalgebra, with usual complex conjugation $\overline{f(x,\xi)}:=\sum_I \overline{f_I(x)}\xi^I$ as superinvolution, represented by multiplication operators on the Hilbert superspace $L^2(\gR^{m|n})$. The choice of $J$ as in Example \ref{ex-l2} leads to the C*-norm $\norm f\norm:=\sum_I \norm f_I\norm_\infty$. But the involution * (attached to the scalar product $\langle-,J(-)\rangle$) does not stabilize $L^\infty(\gR^{m|n})$ so that it is not a C*-algebra. Noncommutative examples of C*-superalgebras will be given in Section \ref{sec-dq}.
\end{example}

Symmetries of noncommutative spaces are quantum groups \cite{Woronowicz:1987,Majid:1995}. In the present graded context, we call the symmetries of noncommutative superspaces quantum supergroups. This notion already exists in the literature at the algebraic level (see e.g. \cite{Majid:1995}) but in the present program, we want to extend this notion at the topological level. Natural objects to consider would be C*-quantum supergroups (with a C*-superalgebra structure) and it is currently under study. A structure going in this direction is that of Fr\'echet quantum supergroups.

\begin{definition}[\cite{deGoursac:2014kv}]
A {\defin Fr\'echet quantum supergroup} is a $\gZ_2$-graded Fr\'echet space with a graded topological tensor product and a continuous homogeneous of degree 0 Hopf algebra structure.

A {\defin representation} of a Fr\'echet quantum supergroup $H$ is a Fr\'echet $\gZ_2$-graded $H$-comodule algebra with continuous homogeneous of degree 0 coaction.
\end{definition}

\begin{example}
The simplest and well-known example corresponds to the Abelian supergroup $\gR^{m|n}$. It can be seen in particular as a Fr\'echet quantum supergroup but commutative by considering $H=\caC^\infty(\gR^{m|n})$, which is a Fr\'echet-Hopf algebra for the usual pointwise product, unit and for the usual coproduct, counit and antipode given by, for any $z_i\in\gR^{m|n}$,
\begin{equation*}
\Delta(f)(z_1,z_2)=f(z_1+z_2),\qquad \eps(f)=f(0),\qquad S(f)(z)=f(-z).
\end{equation*}
We will see noncommutative examples of Fr\'echet quantum supergroups in Section \ref{sec-dq}.
\end{example}

\section{Application to harmonic analysis of Lie supergroups}

Let us present the Heisenberg supergroup. We endow the superspace $\gR^{2m|2n}$ with an even symplectic superform, associated to the matrix
\begin{equation*}
\omega=\begin{pmatrix} 0 & \gone_m & 0 & 0 \\ -\gone_m & 0 & 0 & 0\\ 0 & 0 & \gone_{r+2s} & 0\\ 0 & 0 & 0 & -\gone_r\end{pmatrix}
\end{equation*}
of size $(2m+2n)$ in the canonical basis, with $2n=2r+2s$. the odd part of this symplectic form corresponds to a symmetric form and one has to choose its signature $(r+2s,r)$. The Heisenberg supergroup is the supergroup $G=\gR^{2m|2n}\times \gR^{1|0}Z$ (with $Z$ an even generator) defined as usual by its Lie superalgebra relations
\begin{equation*}
\forall a,a'\in\gR^{2m|2n},\ \forall t,t'\in\gR^{1|0}\quad:\quad [a+tZ,a'+t'Z]=\omega(a,a')Z,
\end{equation*}
where $[-,-]$ denotes the graded Lie bracket.

\begin{proposition}[\cite{Bieliavsky:2010su,deGoursac:2014svn}]
\label{prop-schro}
Kirillov's Orbit Method applied to the Heisenberg supergroup $G$ yields a Hilbert superspace $\ehH_\theta:=L^2(\gR^{m|r})\otimes Hol(\gC^{0|s})$, with elements $\varphi(q,\xi,\zeta)=\sum_{I,J}\varphi_{IJ}(q)\xi^I\zeta^J$ ($q\in\gR^m$, $\xi\in\gR^{0|r}$, $\zeta\in\gC^{0|s}$) and superhermitian inner product
\begin{equation}
\langle \varphi,\psi\rangle=(2i)^s\int \dd q\dd\xi\dd\zeta\dd\overline\zeta\ \overline{\varphi(q,\xi,\zeta)}\psi(q,\xi,\zeta) e^{\frac{i}{\theta}\zeta\overline{\zeta}}.\label{eq-innerprod}
\end{equation}
Moreover, this method also produces a representation $U_\theta$ of $G$ on $\ehH_\theta$
\begin{equation*}
U_\theta(g)\varphi(q_0,\xi_0,\zeta_0)=e^{\frac{2i}{\theta}\big(t+(\frac12 q-q_0)p+(\frac12\xi-\xi_0)\eta+\frac12(\frac12\zeta-\zeta_0)\overline{\zeta}\big)}\varphi(q_0-q,\xi_0-\xi,\zeta_0-\zeta)
\end{equation*}
for $g=(q,p,\xi,\eta,\zeta,\overline\zeta,t)\in G$ with $t\in\gR$, $p\in\gR^m$, $\eta\in\gR^{0|r}$ in the real polarization; $U_\theta$ generalizing the usual Schr\"odinger representation (corresponding to the case $r=s=0$). This representation $U_\theta$ is {\defin superunitary}, namely
\begin{equation*}
\forall g\in G,\ \forall \varphi,\psi\in\ehH_\theta\quad:\quad \langle U_\theta(g)\varphi,U_\theta(g)\psi\rangle=\langle \varphi,\psi\rangle.
\end{equation*}
\end{proposition}

In case of positive signature $r=0$, $n=s$, then the standard inner product \eqref{eq-innerprod2} corresponds to \eqref{eq-innerprod} for which a fundamental symmetry is given by $J(\varphi)=\sum_I (-i)^{|I|}\varphi_I \zeta^I$. However, in case of mixed signature $r\neq 0$, the standard inner product \eqref{eq-innerprod2} does not correspond to \eqref{eq-innerprod} for which any fundamental symmetry has to be more general than a multiplication by a power of $i$. By using superunitarity induced by \eqref{eq-innerprod} and the notion of Hilbert superspace of Definition \ref{def-hs}, one obtains an extended Stone--von Neumann theorem in any signature \cite{deGoursac:2014svn}, identifying the superunitary dual of $G$ with irreducible representations of the Clifford algebra. On the contrary, the standard notion of Hilbert superspace associated to \eqref{eq-innerprod2} leads to an empty dual superunitary in mixed signature \cite[Theorem 5.8]{Salmasian:2010}, which is not satisfactory.
\begin{theorem}[\cite{deGoursac:2014svn}]
Any superunitary representation of $G$ with fixed character $e^{\frac{2i}{\theta}t}$ decomposes as an orthogonal direct sum of the irreducible superunitary representations $U_\theta$.
\end{theorem}
So we see that this structure of Hilbert superspace is more adapted than the standard one \cite{Deligne:1999su,Carmeli:2006} of \eqref{eq-innerprod2} to the analytic point of view of unitary representations of Lie supergroups.

\section{Application to non-formal deformation quantization}
\label{sec-dq}

Contrary to the formal case extensively studied and classified in Refs. \cite{Bayen:1978,Omori:1991,Lecomte:1992,Fedosov:1994,Kontsevich:2003}, there are only few available examples of non-formal deformation quantization of Lie groups in the smooth non-graded setting. Rieffel \cite{Rieffel:1993} built the deformation of Abelian Lie groups and the associated universal deformation formula (UDF). This was also recently extended to (non-Abelian) K\"ahlerian Lie groups \cite{Bieliavsky:2002,Bieliavsky:2010kg} and to the case of $SL(2,\gR)$ \cite{Bieliavsky:2008mv} and of $SU(1,n)$ \cite{Korvers:2014}. Star-exponentials \cite{Cahen:1985,Arnal:1990} associated to these deformations were computed in the non-formal setting in \cite{Bieliavsky:2013sk,Bieliavsky:2013sl2} and such deformations were linked to Hilbert algebras and multipliers in \cite{deGoursac:2014mu}. Non-formal deformation quantization was also extended to the complex case \cite{Omori:2000,Garay:2013gya} and to Abelian $p$-adic groups \cite{Gayral:2014ad}.

Let us examine the deformation quantization of the Heisenberg supergroup $G=\gR^{2m|n}\times \gR^{1|0}Z$. By using the superunitary representation of $G$ from Proposition \ref{prop-schro} as well as the symmetric structure of its coadjoint orbit $\gR^{2m|n}$, we have constructed a non-formal deformation quantization of $\gR^{2m|n}$ invariant under the symmetry $G$.
\begin{theorem}[\cite{Bieliavsky:2010su}]
The star-product defined by
\begin{equation}
(f_1\star_\theta f_2)(z)=\frac{1}{\pi^{2m}\theta^{2m-n}}\int\dd z_1\dd z_2\ K_\theta(z_1,z_2) f_1(z+z_1)f_2(z+z_2),\label{eq-prod-moy}
\end{equation}
for $z=(x,\xi)\in\gR^{2m|n}$, with supergroup cocycle $K_\theta(z_1,z_2)=e^{\frac{-2i}{\theta}\omega(z_1,z_2)}$, on the space
\begin{equation*}
\caB^1(\gR^{2m|n}):=\big\{ f\in\caC^\infty(\gR^{2m|n}),\ \forall \alpha\in\gN^{2m},\ \sup_{x\in\gR^{2m}}\sum_I |\partial_x^\alpha f_I(x)|<+\infty\big\}\quad\simeq\quad\caB(\gR^{2m})\otimes\bigwedge\gR^n,
\end{equation*}
is associative, $G$-invariant and satisfies the tracial property
\begin{equation*}
\int \dd z\ (f_1\star_\theta f_2)(z)=\int \dd z\ f_1(z)f_2(z).
\end{equation*}
Moreover, there is a Weyl-type quantization map $\Omega_\theta:\caB^1(\gR^{2m|n})\to\caB(\ehH_\theta)$, compatible with the star-product and $G$-equivariant, i.e.
\begin{equation*}
\Omega_\theta(f_1\star_\theta f_2)=\Omega_\theta(f_1)\Omega_\theta(f_2),\qquad \Omega_\theta(L^*_{g^{-1}}f)=U_\theta(g)\Omega_\theta(f) U_\theta(g)^{-1}.
\end{equation*}
\end{theorem}
The space $\caB^1(\gR^{2m|n})$, endowed with the star-product, is a Fr\'echet algebra isomorphic to the tensor product of the Moyal--Weyl algebra $\caB(\gR^{2m})$ and of the Clifford algebra $Cl(n,\gC)$. Note that the space $\caB^1(\gR^{2m|n})$ is part of a pseudodifferential calculus adapted to the superspace $\gR^{2m|n}$ for which an oscillatory integral is also defined.

An analogous construction was also performed \cite{deGoursac:2014sp} with a different cocycle for the {\defin Poincar\'e supergroup} in $2|2$ dimensions, i.e. the supergroup defined by even generators $H,E,F$, odd generators $\Gamma,\Xi$ and Lie superalgebra relations
\begin{equation*}
[H,E]=2E,\quad [H,F]=-2F,\quad [H,\Gamma]=\Gamma,\quad [H,\Xi]=-\Xi,\quad [\Gamma,\Gamma]=-2E,\quad [\Xi,\Xi]=-2F.
\end{equation*}

For the Heisenberg supergroup, a universal deformation formula was also built, namely the use of $\star_\theta$, can also deform any algebra on which the supergroup $\gR^{2m|n}$ is acting. Let us describe it. We consider an action $\rho$ of the supergroup $\gR^{2m|n}$ on a Fr\'echet algebra $(\algA,|\fois|_j)$ satisfying the conditions:
\begin{itemize}
\item $\forall z_1,z_2\in\gR^{2m|n}$, $\rho_{z_1+z_2}=\rho_{z_1}\rho_{z_2}$ and $\rho_0=\text{id}$.
\item $\forall z\in\gR^{2m|n}$, $\forall a,b\in\algA$, $\rho_z(ab)=\rho_z(a)\rho_z(b)$, and $\rho_z$ is linear.
\item By writing $z=(x,\xi)\in\gR^{2m|n}$, we can expand the action as: $\rho_{(x,\xi)}(a)=\sum_I \rho_x(a)_I\xi^I$; $\forall a\in\algA$, $\forall I$, $x\mapsto \rho_x(a)_I$ is $\algA$-valued and continuous.
\item The action is subisometric, i.e. there exists a constant $C>0$ such that $\forall a\in\algA$, $\forall I$, $\forall j$, $\exists k$, $\forall x\in\gR^{2m}$, $|\rho_x(a)_I|_j\leq C|a|_k$.
\end{itemize}
The star-product \eqref{eq-prod-moy} can be directly extended to $\algA$-valued superfunctions $\caB^1_\algA(\gR^{2m|n})$. Note that this space is Fr\'echet for the seminorms $|f|_{j,\alpha} =\sup_{x\in\gR^{2m}}\{\sum_I|\partial_x^\alpha f_I(x)|_j\}$. We recall that the set of {\defin smooth vectors} of $\algA$ for the action $\rho$ is defined as
\begin{equation*}
\algA^\infty=\{a\in\algA,\quad \rho^a:=z\mapsto \rho_z(a)\text{ is smooth on }\gR^{2m|n}\}.
\end{equation*}
This set $\algA^\infty$ is dense in $\algA$, and for any $a\in\algA^\infty$, the map $\rho^a$ lies in $\caB^1_{\algA^\infty}(\gR^{2m|n})$, so that we can form the star-product of $\rho^a$ and $\rho^b$, for $a$ and $b$ smooth vectors.
\begin{theorem}[\cite{Bieliavsky:2010su}]
\label{thm-udf}
The expression
\begin{equation*}
a\star_\theta b:=(\rho^a\star_\theta\rho^b)(0),
\end{equation*}
for $a,b\in\algA^\infty$, defines an associative product on $\algA^\infty$. If $\algA$ is a C*-superalgebra and $\rho$ is compatible with the degree and the superinvolution, then $(\algA^\infty,\star_\theta)$ can be completed to another C*-superalgebra denoted by $\algA_\theta$.
\end{theorem}
Therefore, the category of C*-superalgebras is stable by this deformation quantization, which stresses the interest of this structure for noncommutative supergeometry. Note that if $\algA$ is a C*-algebra on which $\gR^{2m|n}$ acts, the deformation $(\algA^\infty,\star_\theta)$ cannot be completed into a C*-algebra in general. The {\defin universal deformation formula} (UDF) of Theorem \ref{thm-udf} is a generalization of Drinfeld twists \cite{Drinfeld:1989,Giaquinto:1998} to the non-formal setting.

This UDF thus produces a large class of (deformed) noncommutative superspaces. For example, as soon as $\gR^{2m|n}$ acts on a trivial compact supermanifold, one can deform its algebra of continuous superfunctions to a noncommutative C*-superalgebra \cite{Bieliavsky:2010su}. Noncommutative supertori appear as particular cases of this deformation.
\begin{definition}[\cite{Bieliavsky:2010su}]
The {\defin noncommutative supertorus} of dimension $2m|n$ and signature $(p,q)$ (with $n=p+q$) is the C*-superalgebra generated by the even elements $U_j$, $V_j$ ($1\leq j\leq m$) and by the odd elements $\Gamma_k$, $\Xi_\ell$ ($1\leq k\leq p$, $1\leq \ell\leq q$) satisfying the relations
\begin{align*}
& U_j^\dag U_j=U_jU_j^\dag=V_j^\dag V_j=V_jV_j^\dag=\gone,\qquad \Gamma_k^\dag=\Gamma_k,\qquad \Xi_\ell^\dag=\Xi_\ell,\\
& U_jV_j= e^{2i\pi\theta}V_jU_j,\qquad (\Gamma_k)^2=i\theta,\qquad (\Xi_\ell)^2=-i\theta,
\end{align*}
the other relations being trivial. It coincides with a deformation of the usual supertorus of dimension $2m|n$, with deformation parameter $\theta$.
\end{definition}

\medskip

By using the UDF of Theorem \ref{thm-udf}, we can also construct interesting examples of Fr\'echet quantum supergroups. Inspired by \cite{Rieffel:1993qg} in the non-graded case, we consider an action $\pi:\gR^{1|0}\to Sp(\gR^{2m|n},\omega)$ and $G':=\gR^{1|0}\ltimes_\pi\gR^{2m|n}$. $G'$ is a solvable Lie supergroup.
\begin{theorem}[\cite{deGoursac:2014kv}]
The space $H:=\caC^\infty(\gR^{1|0})\hat\otimes\caB^1(\gR^{2m|n})$, endowed with the usual Hopf algebra structure of $G'$ except for the product, which is given by the deformation
\begin{equation*}
f_1\star f_2(t,z):=\frac{1}{\pi^{2m}t^{2m-n}}\int\dd z_1\dd z_2\ f_1(t,z+z_1)f_2(t,z+z_2) K_t(z_1,z_2),
\end{equation*}
for any $t\in\gR^{1|0}\simeq\gR$ and $z\in\gR^{2m|n}$, is a Fr\'echet quantum supergroup. An analogue of the {\defin multiplicative unitary} operator $W\in\caL(H\hat\otimes H)$ can be defined for this quantum supergroup by
\begin{equation*}
W(f_1\otimes f_2)=(\Delta f_1)\star (1\otimes f_2).
\end{equation*}
It satisfies the pentagonal equation $W_{12}W_{13}W_{23}=W_{23}W_{12}$ as well as the superunitarity condition with respect to a Hilbert superspace structure of $L^2(G'\times G')$.
\end{theorem}
A similar construction of Fr\'echet quantum supergroup was performed \cite{deGoursac:2014kv} for connected Lie supergroups with supertoral subgroups.

\section{Application to QFT}

We give here an application where a NCSG is useful to understand a renormalizable noncommutative quantum field theory. Let us consider the Moyal space \cite{GraciaBondia:1987kw,Cagnache:2009ik} given by the deformation quantization of $\gR^{2m}$. It turns out that the $\phi^{\star_\theta 4}$ theory on the Moyal space is not renormalizable because of a new divergence called UV-IR mixing \cite{Minwalla:1999px}, generic to noncommutative spaces \cite{Gayral:2005af}. The first solution of this problem was found by H. Grosse and R. Wulkenhaar by adding a {\defin harmonic term} to the standard {\defin $\phi^{\star_\theta 4}$ action functional}
\begin{equation}
S_{GW}[\phi]=\int \Big(\frac 12(\partial_\mu\phi)^2+\frac{2\Omega^2}{\theta^2} x^2\phi^2+\frac{M^2}{2}\phi^2+\lambda\, \phi^{\star_\theta 4}\Big)\dd x,\label{eq-actharm}
\end{equation}
which becomes a renormalizable model ($m=1,2$) to all orders \cite{Grosse:2003nw,Grosse:2004yu,Rivasseau:2005bh} for $\Omega\neq 0$. This model possesses new interesting properties concerning its vacuum configurations \cite{deGoursac:2007uv}, its Connes--Kreimer algebra \cite{Tanasa:2007xa}, its symmetries \cite{deGoursac:2009fm,Hounkonnou:2009qt}, its commutative limit \cite{deGoursac:2012ki}, its beta function \cite{Disertori:2006nq} and finally its solvability for $\theta\to\infty$ \cite{Grosse:2014lxa}. Note that there exists now another renormalizable real scalar theory on the Moyal space \cite{Gurau:2008vd,Tanasa:2010fk}. A gauge theory associated to \eqref{eq-actharm} was constructed in \cite{deGoursac:2007gq,Grosse:2007dm}. It exhibits similar features \cite{deGoursac:2008rb,Cagnache:2008tz} and is a candidate to renormalizability. See \cite{Blaschke:2007vc,Blaschke:2009aw,Blaschke:2013gha} for the study of the BRST symmetry.
\medskip

We want to give an interpretation of the harmonic term in \eqref{eq-actharm}. It turns out that the action \eqref{eq-actharm} can be reformulated in terms of the commutator and anticommutator
\begin{equation*}
S_{GW}[\phi]=\int \Big(\frac 12\big|[\frac{i}{2}x_\mu,\phi]_{\star_\theta}\big|^2+ \frac{\Omega^2}{2} \big|\{\frac{i}{2}x_\mu,\phi\}_{\star_\theta}\big|^2+\frac{M^2}{2}\phi^2+\lambda\, \phi^{\star_\theta 4}\Big)\dd x.
\end{equation*}
This is a sign of a graded symmetry. Actually, this model can be interpreted in the setting of deformation quantization of the superspace $\gR^{2m|1}$. We consider the natural differential calculus \cite{deGoursac:2008bd} for this graded algebra associated to the graded derivations $\mathfrak{d}\in\big\{[\frac{i}{2}x_\mu,-]_{\star_\theta},\ [\frac{i}{2}x_\mu \xi,-]_{\star_\theta}\big\}_{ \mu\in\{1,..,2m\}}$.
\begin{proposition}[\cite{deGoursac:2010zb,Bieliavsky:2010su}]
The standard $\Phi^{\star_\theta 4}$ action of the deformation quantization of $\gR^{2m|1}$, after integration with respect to the odd variable $\xi$ and identification of the even and odd fields $\phi_1(x)=b\phi_0(x)$ ($b\in\gR$) in $\Phi(x,\xi)=\phi_0(x)+\phi_1(x)\xi$, gives exactly rise to the renormalizable Grosse--Wulkenhaar action \eqref{eq-actharm}:
\begin{equation*}
\text{tr}\Big(\frac12\sum_{\mathfrak{d}}|\mathfrak{d}(\Phi)|^2+\frac{M^2}{2}\Phi^2+\Lambda\Phi^{\star_\theta 4}\Big)=S_{GW}[\phi_0]
\end{equation*}
with $\Omega^2=\frac{b^4\theta^2}{16}$ and $\lambda=\Lambda(1+\frac{b^4\theta^2}{16})$.

In this setting, the Langmann--Szabo duality \cite{Langmann:2002cc} corresponds to the switch of grading and the gauge action of \cite{deGoursac:2007gq} can be directly derived as the Yang--Mills action of the deformation of $\gR^{2m|1}$ with differential calculus associated to the graded derivations $\mathfrak d$.
\end{proposition}
This shows that the NCSG of the deformation of $\gR^{2m|1}$ is likely to be related to the renormalizability of the noncommutative quantum field theory \eqref{eq-actharm}. This interpretation could be useful to obtain renormalizable quantum field theories on other noncommutative spaces.

\bibliographystyle{utcaps}
\bibliography{biblio-perso,biblio-recents,biblio-these}

\end{document}